\documentclass[11pt]{article}
\usepackage{amscd}
\usepackage{amsfonts}
\usepackage{amsmath}
\usepackage{amssymb}
\usepackage{amsthm}
\usepackage{bbm}
\usepackage{CJK}
\usepackage{fancyhdr}
\usepackage{graphicx}
\usepackage{hyperref}
\usepackage{indentfirst}
\usepackage{latexsym}
\usepackage{mathrsfs}
\usepackage{xypic}
\usepackage{amsmath,amssymb,amsfonts,amsthm,fancyhdr}
\usepackage{epsfig,graphicx,picins,picinpar,subfigure}
\usepackage{pstricks}
\usepackage{amssymb}
\usepackage{xypic}
\usepackage[compress]{cite}

\newtheorem{theorem}{Theorem}[section]
\newtheorem{lemma}[theorem]{Lemma}
\newtheorem{definition}[theorem]{Definition}
\newtheorem{proposition}[theorem]{Proposition}
\newtheorem{example}[theorem]{Example}

\usepackage[top=1in,bottom=1in,left=1.25in,right=1.25in]{geometry}
\def\<{\langle}
\def\>{\rangle}

\def\c{\cdot}

\date{}
\begin{document}
\renewcommand{\baselinestretch}{1.2}
\renewcommand{\arraystretch}{1.0}
\title{\bf Embedding tensors on 3-Hom-Lie algebras}
\author{{\bf Wen Teng$^{1}$,   Jiulin Jin$^{2}$, Yu Zhang$^{1}$}\\
{\small 1. School of Mathematics and Statistics, Guizhou University of Finance and Economics} \\
{\small  Guiyang  550025, P. R. of China}\\
  {\small Email: tengwen@mail.gufe.edu.cn(Wen Teng)} \\
{\small 2.   College of Mathematics and Information Science,  Guiyang University}\\
{\small Guizhou  $550005$, P. R. of China}\\}
 \maketitle
\begin{center}
\begin{minipage}{13.cm}

{\bf Abstract}  In this paper,      we introduce the notion of   embedding tensor on   3-Hom-Lie algebras and naturally induce  3-Hom-Leibniz algebras. Moreover,   the cohomology theory
of embedding tensors on 3-Hom-Lie algebras is defined. As an application, we show that if two linear deformations of an embedding tensor on a 3-Hom-Lie algebra are equivalent, then their infinitesimals belong to
the same cohomology class in the first cohomology group.
 \smallskip

 \smallskip

 \end{minipage}
 \end{center}
 \normalsize\vskip0.5cm

\section{Introduction}
\def\theequation{\arabic{section}. \arabic{equation}}
\setcounter{equation} {0}

The concept of embedded tensors provides a useful tool in the construction of supergravity theories\cite{Nicolai} and higher gauge theories\cite{Bergshoeff}.
In mathematics, embedding tensor is called   average operator. Aguiar \cite{Aguiar} studies the average operator on the associative algebra    and Lie algebra.
Further,   deformation and cohomology theory of embedding tensors of associative algebra, Lie algebra and 3-Lie algebra were given in \cite{Das1,Sheng,Hu}. Recently,  Das and Makhlouf\cite{Das} introduced the embedding tensor on Hom-Lie algebra, and studied the related properties.

The aim of this paper is to extend the concept of embedded tensors of 3-Lie algebras to Hom-type algebras.
Hom-Lie algebras were introduced by Hartwig, Larsson and Silvestrov \cite{Hartwig}
in the study of q-deformations  of the Witt and Virasoro algebra.
  In the last fifteen years,  Hom-type algebras have attracted extensive attention from scholars (see \cite{Benayadi,Chen2012,Das,Makhlouf,Makhlouf2010,Mabrouk,Mishraa, Song}). On the other hand, Filippov\cite{Filippov}  introduced 3-Lie algebra and more general $n$-Lie algebra, which can be regarded as a generalization of Lie algebra to   higher  algebra.  In \cite{ Ataguema}, $n$-Hom-Lie algebras and various generalizations of n-ary algebras are considered, and the representation and cohomology of $n$-Hom-Lie algebras are first studied in \cite{Ammar}.

  This paper is organized as follows.
In Section 2, we  recall the definitions of  3-Hom-Lie algebras. Then we introduce the notion of an embedding tensor on a 3-Hom-Lie algebra, which naturally induces a 3-Hom-Leibniz algebra.
In Section 3, we  introduce  the representation and cohomology theory of embedding tensor on a 3-Hom-Lie algebra.
In Section 4, we study linear deformations of embedding tensor on a 3-Hom-Lie algebra, we show that if two linear deformations of an embedding tensor on a 3-Hom Lie algebra are equivalent, then their infinitesimals  belong
to
the same cohomology class in the first cohomology group.

 In this paper, all vector spaces are considered over a field $ \mathbb{K}$ of characteristic 0.

\section{Embedding tensors on 3-Hom-Lie algebras}
\def\theequation{\arabic{section}.\arabic{equation}}
\setcounter{equation} {0}

In this section, we recall some basic definitions of 3-Hom-Lie algebras,  Hom-Leibniz algebras and 3-Hom-Leibniz algebras. Then we introduce embedding tensors on 3-Hom-Lie algebras. We show that an embedding
tensor naturally gives rise to a 3-Hom-Leibniz algebra structure.  Finally, we provide some examples of embedding tensors on
3-Hom-Lie algebras.

\begin{definition}(\cite{Ataguema})  A 3-Hom-Lie algebra is a triple $(L, [\c,\c,\c],  \alpha)$ consisting of a vector space $L$,  a trilinear
skew-symmetric mapping $[\c,\c,\c]: L\times L \times L \rightarrow L$,  and a linear map $\alpha: L \rightarrow L$   satisfying $\alpha([x,y,z])=[\alpha(x),\alpha(y),\alpha(z)]$ and the Hom-Filippov-Jacobi
identity:
\begin{align}
&[\alpha(a),\alpha(b),[x,y,z]]=[[a,b,x], \alpha(y), \alpha(z)]+ [\alpha(x),[a,b,y],\alpha(z)]+[\alpha(x),\alpha(y),[a,b,z]]\label{2.1},
\end{align}
for any  $a,b, x, y, z\in L$.
Furthermore, if $\alpha: L \rightarrow L$  is a vector space automorphism of $L$, then the 3-Hom-Lie algebra $(L, [\c,\c,\c],  \alpha)$ is called a regular 3-Hom-Lie algebra.
\end{definition}

A homomorphism between two 3-Hom-Lie algebras  $(L, [\cdot, \cdot,\c], \alpha)$ and $(L', [\cdot, \cdot,\c]',  \alpha')$ is a linear map $\psi: L\rightarrow L'$ satisfying $\psi\circ\alpha=\alpha'\circ\psi$ and
\begin{eqnarray*}
&&\psi([x, y,z])=[\psi(x), \psi(y), \psi(z)]', ~~\forall ~x,y,z\in L.
\end{eqnarray*}
In particular, if $\psi$ is nondegenerate, then $\psi$ is called an isomorphism from $L$ to $L'$.

\begin{definition} (\cite{Makhlouf2010})
A   Hom-Leibniz algebra is a vector space $L$ together with a
bracket operation $[\c,\c]: L\times L \rightarrow L$ and a linear map $\alpha: L \rightarrow L$   satisfying $\alpha([x,y])=[\alpha(x),\alpha(y)]$   such that
\begin{eqnarray*}
&&[\alpha(x),[y,z]]= [[x,y], \alpha(z)]+ [\alpha(y),[x,y]],
\end{eqnarray*}
for any  $x,y,z\in L$.

\end{definition}

Let $(L, [\c,\c,\c],  \alpha)$ be a 3-Hom-Lie algebra, then the elements in $\wedge^2 L$ are called fundamental objects of the 3-Hom-Lie algebra $(L, [\c,\c,\c],  \alpha)$. If there is a bilinear operation $[\c,\c]$ on $\wedge^2 L$, which is given by
\begin{eqnarray*}
&&[X,Y]'=[x_1,x_2,y_1]\wedge \alpha(y_2)+\alpha(y_1)\wedge[x_1,x_2,y_2], \forall X=x_1\wedge x_2, Y=y_1\wedge y_2\in  \wedge^2 L
\end{eqnarray*}
and a linear map $\tilde{\alpha}$ on $\wedge^2 L$ is defined by $\tilde{\alpha}(X)=\alpha(x_1)\wedge \alpha(x_2)$.
Clearly, $(\wedge^2 L, [\c,\c]', \tilde{\alpha})$ is a Hom-Leibniz algebra\cite{Daletskii}.

\begin{definition} (\cite{Mabrouk})
A representation of a 3-Hom-Lie algebra $(L, [\c,\c,\c],  \alpha)$ on a vector space $V$ with respect
to $\beta\in \mathrm{End}(V)$ is a skew-symmetric linear map $\rho:\wedge^2 L\rightarrow \mathrm{End}(V)$ such that
\begin{align}
&\rho(\alpha(x),\alpha(y))\circ\beta=\beta\circ\rho(x,y),\label{2.2}\\
&\rho(\alpha(x),\alpha(y))\rho(a,b)-\rho(\alpha(a),\alpha(b))\rho(x,y)=(\rho([x,y,a],\alpha(b))-\rho([x,y,b],\alpha(a)))\circ\beta,\label{2.3}\\
&\rho([x,y,a],\alpha(b))\circ\beta-\rho(\alpha(y),\alpha(a))\rho(x,b)=\rho(\alpha(a),\alpha(x))\rho(y,b)+\rho(\alpha(x),\alpha(y))\rho(a,b),\label{2.4}
\end{align}
for any  $x, y, a, b\in L$.
Furthermore, if $\beta: V\rightarrow V$  is a vector space automorphism of $V$, then  $(V; \rho,  \beta)$ is called a regular representation of  $(L, [\c,\c,\c],  \alpha)$.
\end{definition}

It follows from the above definition that any 3-Hom-Lie algebra  $(L, [\cdot, \cdot,\c], \alpha)$  can be regarded as a representation
of itself, where $\rho=\mathrm{ad}:\wedge^2 L\rightarrow \mathrm{End}(L)$ is given by $\mathrm{ad}(x,y)(z):=[x,y,z]$, for $x,y,z\in L$. This is called the  adjoint
representation.

\begin{definition} (\cite{Makhlouf}) A 3-Hom-Leibniz algebra is a triple $(\mathcal{L}, [\c,\c,\c]_{\mathcal{L}},  \alpha_{\mathcal{L}})$ consisting of a vector space $\mathcal{L}$,  a trilinear
  mapping $[\c,\c,\c]_{\mathcal{L}}: \mathcal{L}\times \mathcal{L} \times \mathcal{L} \rightarrow \mathcal{L}$,  and  a linear map $\alpha_{\mathcal{L}}: \mathcal{L} \rightarrow \mathcal{L}$   satisfying $\alpha_{\mathcal{L}}([x,y,z]_{\mathcal{L}})=[\alpha_{\mathcal{L}}(x),\alpha_{\mathcal{L}}(y),\alpha_{\mathcal{L}}(z)]_{\mathcal{L}}$ such that
\begin{align}
 &[\alpha_{\mathcal{L}}(a),\alpha_{\mathcal{L}}(b),[x,y,z]_{\mathcal{L}}]_{\mathcal{L}}\nonumber\\
 &= [[a,b,x]_{\mathcal{L}}, \alpha_{\mathcal{L}}(y), \alpha_{\mathcal{L}}(z)]_{\mathcal{L}}+ [\alpha_{\mathcal{L}}(x),[a,b,y]_{\mathcal{L}},\alpha_{\mathcal{L}}(z)]_{\mathcal{L}}+[\alpha_{\mathcal{L}}(x),\alpha_{\mathcal{L}}(y),[a,b,z]_{\mathcal{L}}]_{\mathcal{L}},\label{2.5}
\end{align}
for any  $a,b, x, y, z\in L$.
\end{definition}

\begin{proposition}
Let $(L, [\c,\c,\c],\alpha)$ be a 3-Hom-Lie algebra,  $V$ be a  vector space,   $\beta\in \mathrm{End}(V)$ and  $\rho: \wedge^2 L\rightarrow \mathrm{End}(V)$ be  a skew-symmetric linear map. Then $(V; \rho, \beta)$ is a representation of 3-Hom-Lie algebra $L$ if and only if $L \oplus V$ is a 3-Hom-Leibniz algebra under the following maps:
\begin{eqnarray*}
&&(\alpha\oplus\beta)(x+u):=\alpha(x)+\beta(u),\\
&&[x+u, y+v, z+w]_{\rho}:=[x, y, z]+\rho(x, y)(w),
\end{eqnarray*}
for any    $x, y, z\in L$ and $u, v, w\in V$.  $(L \oplus V, [\c,\c,\c]_{\rho},\alpha\oplus\beta)$  is called
the hemisemidirect product 3-Hom-Leibniz algebra, and denoted by $L\ltimes_{\rho}V$.
\end{proposition}

{\bf Proof.}  For all $x,y,z,a,b\in L, u,v,w,s,t\in V,$ by Eqs. \eqref{2.1},  \eqref{2.2}  and \eqref{2.3} , we have
\begin{align*}
&(\alpha\oplus\beta)([x+u,y+v,z+w]_{\rho})=(\alpha\oplus\beta)([x, y, z]+\rho(x, y)(w))\\
=&\alpha([x, y, z])+\beta(\rho(x, y)(w))=[\alpha(x), \alpha(y), \alpha(z)]+\rho(\alpha(x), \alpha(y))\beta(w)\\
=&[\alpha(x)+\beta(u), \alpha(y)+\beta(v), \alpha(z)+\beta(w)]_{\rho}\\
=&[(\alpha\oplus\beta)(x+u), (\alpha\oplus\beta)(y+v), (\alpha\oplus\beta)(z+w)]_{\rho},\\
&[[a+s,b+t,x+u]_{\rho}, (\alpha\oplus\beta)(y+v), (\alpha\oplus\beta)(z+w)]_{\rho}\\
&+ [(\alpha\oplus\beta)(x+u),[a+s,b+t,y+v]_{\rho},(\alpha\oplus\beta)(z+w)]_{\rho}\\
&+[(\alpha\oplus\beta)(x+u),(\alpha\oplus\beta)(y+v),[a+s,b+t,z+w]_{\rho}]_{\rho}\\
&-[(\alpha\oplus\beta)(a+s),(\alpha\oplus\beta)(b+t),[x+u,y+v,z+w]_{\rho}]_{\rho}\\
=&[[a, b, x]+\rho(a, b)(u), \alpha(y)+\beta(v),\alpha(z)+\beta(w)]_{\rho}\\
& +[\alpha(x)+\beta(u),[a, b, y]+\rho(a, b)(v),\alpha(z)+\beta(w)]_{\rho}\\
&+[\alpha(x)+\beta(u),\alpha(y)+\beta(v),[a, b, z]+\rho(a, b)(w)]_{\rho}\\
&-[\alpha(a)+\beta(s),\alpha(b)+\beta(t),[x, y, z]+\rho(x, y)(w)]_{\rho}\\
=&[[a, b, x], \alpha(y),\alpha(z)]+\rho([a, b, x], \alpha(y))\beta(w)+[\alpha(x),[a, b, y],\alpha(z)]+\rho(\alpha(x),[a, b, y])\beta(w)\\
&+[\alpha(x),\alpha(y),[a, b, z]]+\rho(\alpha(x),\alpha(y))\rho(a, b)(w)-[\alpha(a),\alpha(b),[x, y, z]]-\rho(\alpha(a),\alpha(b))\rho(x, y)(w)\\
=&0.
\end{align*}
Thus,   $(L \oplus V, [\c,\c,\c]_{\rho},\alpha\oplus\beta)$  is a 3-Hom-Leibniz algebra.   \hfill $\square$

\begin{definition}  Let $(V;\rho,\beta)$ be a  representation of the  3-Hom-Lie algebra $(L, [\c,\c,\c],\alpha)$.  Then the linear map
$T:V\rightarrow L$ is called an embedding tensor on the 3-Hom-Lie algebra  $(L, [\c,\c,\c],\alpha)$ with respect to
the representation $(V;\rho,\beta)$ if $T$   meets the following Eqs.:
\begin{align}
&T\circ\beta=\alpha\circ T,\label{2.6}\\
&[Tu,Tv,Tw]=T(\rho(Tu,Tv)w),\label{2.7}
\end{align}
for any  $u,v,w\in V$.
\end{definition}

\begin{theorem}
 A linear map $T:V\rightarrow L$ is an embedding tensor on the 3-Hom-Lie algebra  $(L, [\c,\c,\c],\alpha)$ with respect
to the representation $(V;\rho,\beta)$  if and only if the graph $Gr(T)=\{Tu+u~|~u\in V\}$ is a 3-Hom-Leibniz
subalgebra of the hemisemidirect product 3-Hom-Leibniz algebra $L\ltimes_{\rho}V$.
\end{theorem}

{\bf Proof.}  Let $T:V\rightarrow L$ be a linear map. Then for all  $u,v,w\in V$,  we have
\begin{align*}
&(\alpha\oplus\beta)(Tu+u)=\alpha(Tu)+\beta(u),\\
&[Tu+u,Tv+v,Tw+w]_{\rho}=[Tu,Tv,Tw]+\rho(Tu,Tv)w,
\end{align*}
Thus,    the graph $Gr(T)=\{Tu+u~|~u\in V\}$ is a subalgebra of the hemisemidirect product
3-Hom-Leibniz algebra $L\ltimes_{\rho}V$ if and only if $T$  meets Eqs. \eqref{2.6} and \eqref{2.7}, which implies
that $T$ is an embedding tensor on the 3-Hom-Lie algebra $L$ with respect to the representation $(V;\rho,\beta)$ .
   \hfill $\square$

Clearly, the algebraic structure underlying an embedding tensor on the 3-Hom-Lie algebra $(L, [\c,\c,\c],\alpha)$ is a
3-Hom-Leibniz algebra. Thereby, we have the following proposition.

\begin{proposition}
Let $T:V\rightarrow L$ be an embedding tensor on the 3-Hom-Lie algebra $(L, [\c,\c,\c],\alpha)$ with
respect to the representation $(V;\rho,\beta)$. If a linear map $[\c,\c,\c]_T: V\times V\times V\rightarrow V$ is given by
\begin{align}
&[u,v,w]_T=\rho(Tu,Tv)w, \label{2.8}
\end{align}
for any     $u, v, w\in V$. Then $(V,[\c,\c,\c]_T,\beta)$ is a 3-Hom-Leibniz algebra. Moreover, $T$ is a homomorphism from the 3-Hom-Leibniz
algebra $(V,[\c,\c,\c]_T,\beta)$ to the 3-Hom-Lie algebra $(L,[\c,\c,\c],\alpha)$.
\end{proposition}

{\bf Proof.}  For all $u,v,w,s,t\in V,$ by Eqs. \eqref{2.2},  \eqref{2.3},   \eqref{2.6}, \eqref{2.7} and \eqref{2.8}, we have
\begin{align*}
&\beta([u,v,w]_T)=\beta(\rho(Tu,Tv)w)\\
=&\rho(\alpha(Tu),\alpha(Tv))\beta(w)=\rho(T\beta(u),T\beta(v))\beta(w)\\
=&[\beta(u),\beta(v),\beta(w)]_T,\\
&[[s,t,u]_T, \beta(v), \beta(w)]_T+ [\beta(u),[s,t,v]_T,\beta(w)]_T+[\beta(u),\beta(v),[s,t,w]_T]_T\\
=&[\rho(Ts,Tt)u, \beta(v), \beta(w)]_T+ [\beta(u),\rho(Ts,Tt)v,\beta(w)]_T+[\beta(u),\beta(v),\rho(Ts,Tt)w]_T\\
=&\rho(T\rho(Ts,Tt)u, T\beta(v))\beta(w)+ \rho(T\beta(u),T\rho(Ts,Tt)v)\beta(w)+\rho(T\beta(u),T\beta(v))\rho(Ts,Tt)w\\
=&\rho([Ts,Tt,Tu], \alpha (Tv))\beta(w)+ \rho(\alpha (Tu),[Ts,Tt,Tv])\beta(w)+\rho(\alpha (Tu),\alpha (Tv))\rho(Ts,Tt)w\\
=&\rho(\alpha(Ts),\alpha(Tt))\rho(Tu,Tv)w\\
=&\rho(T\beta(s),T\beta(t))\rho(Tu,Tv)w\\
=&[\beta(s), \beta(t), \rho(Tu,Tv)w]_T\\
=&[\beta(s), \beta(t), [u,v,w]_T]_T
\end{align*}
Thus,   $(V,[\c,\c,\c]_T,\beta)$ is a 3-Hom-Leibniz algebra.  By  Eqs.  \eqref{2.6} and \eqref{2.7}, $T$ is a homomorphism from the 3-Hom-Leibniz
algebra $(V,[\c,\c,\c]_T,\beta)$ to the 3-Hom-Lie algebra $(L,[\c,\c,\c],\alpha)$.   \hfill $\square$

\begin{definition}
Let $T$ and $T'$  be two embedding tensors on the 3-Hom-Lie algebra $(L, [\c,\c,\c],\alpha)$ with
respect to the representation $(V;\rho,\beta)$. Then a homomorphism from $T'$
to $T$ consists of a 3-Hom-Lie algebra
homomorphism $\psi_L:L\rightarrow L$ and a linear map $\psi_V:V\rightarrow V$ such that
\begin{align}
 &\beta\circ \psi_V=\psi_V\circ \beta,\label{2.9}\\
 &T\circ \psi_V=\psi_L\circ T',\label{2.10}\\
 &\psi_V(\rho(x,y)u)=\rho(\psi_L(x),\psi_L(y))\psi_V(u).\label{2.11}
\end{align}
In particular, if both $\psi_L$ and $\psi_V$ are invertible, $(\psi_L,\psi_V)$ is called an isomorphism from $T'$
to $T$.
\end{definition}

The association of a 3-Hom-Leibniz algebra from a   embedding tensor enjoys the functorial
property.

\begin{proposition}
Let $T$ and $T'$  be two embedding tensors on the 3-Hom-Lie algebra $(L, [\c,\c,\c],\alpha)$ with
respect to the representation $(V;\rho,\beta)$, and $(\psi_L,\psi_V)$  be a homomorphism from $T'$
to $T$. Then $\psi_V$  is a
homomorphism of 3-Hom-Leibniz algebras from $(V,[\c,\c,\c]_{T'},\beta)$ to $(V,[\c,\c,\c]_T,\beta)$.
\end{proposition}

{\bf Proof.}  For all $u,v,w\in V,$ by Eqs. \eqref{2.8}, \eqref{2.10} and \eqref{2.11} , we have
\begin{align*}
&\psi_V([u,v,w]_{T'})\\
=&\psi_V(\rho(T'u,T'v)w)\\
=&\rho(\psi_L(T'u),\psi_L(T'v))\psi_V(w)\\
=&\rho(T\psi_V(u),T\psi_V(v))\psi_V(w)\\
=&[\psi_V(u),\psi_V(v),\psi_V(w)]_T.
\end{align*}
Using Eq. \eqref{2.9}, we can get $\psi_V$  is a
homomorphism of 3-Hom-Leibniz algebras from $(V,[\c,\c,\c]_{T'},\beta)$ to $(V,[\c,\c,\c]_T,\beta)$.   \hfill $\square$

Next, we present some examples of embedding tensors on 3-Hom-Lie algebras.

\begin{example}
Let $(L, [\c,\c,\c],\alpha)$ be a 3-Hom-Lie algebra. Then the identity map $\mathrm{Id}:L\rightarrow L$ is an embedding tensor
on the 3-Hom-Lie algebra $L$ with respect to the adjoint representation.
\end{example}

\begin{example}
Let $(L, [\c,\c,\c],\alpha)$ be a 3-Hom-Lie algebra. Then a linear map $D:L\rightarrow L$ is said to be a derivation for the
3-Hom-Lie algebra $L$ if $\alpha\circ D=D\circ\alpha$ and $D[x,y,z]=[Dx,y,z]+[x,Dy,z]+[x,y,Dz]$, for all $x,y,z\in L$. If $D^2=0,$   then $D$ is an embedding tensor on $L$ with respect to the adjoint representation.
\end{example}

\begin{example}
Let $(V;\rho,\beta)$ be a  representation of a 3-Hom-Lie algebra$(L, [\c,\c,\c],\alpha)$.  If a  linear map
$f:V\rightarrow L$  satisfies:
\begin{align*}
&\alpha(f(u))=f(\beta(u)),\\
&f(\rho(x,f(u))v)=[x,f(u),f(v)],
\end{align*}
for any $x\in L, u,v\in V.$ Then $f$ is an embedding tensors on the 3-Hom-Lie algebra $(L, [\c,\c,\c],\alpha)$ with
respect to the representation $(V;\rho,\beta)$
\end{example}

\begin{example}
Let $(L, [\c,\c,\c],\alpha)$ be a 3-Hom-Lie algebra. Then   it can be easily checked that $(\oplus^n L; \rho, \oplus^n \alpha) $ is
a representation of the 3-Hom-Lie algebra   $L$, where
$$\rho:\wedge^2 L\rightarrow \mathrm{End}(\oplus^n L), \rho(x,y)((x_1,\cdots,x_n))=([x,y,x_1],\cdots,[x,y,x_n]),$$
for any $x\in L, (x_1,\cdots, x_n)\in \oplus^n L$. Moreover, $T:\oplus^n L\rightarrow L, T(x_1,\cdots, x_n)=x_1+\cdots+x_n$ is an embedding tensor on $L$ with respect to the  representation $(\oplus^n L; \rho, \oplus^n \alpha) $.
\end{example}

\begin{example}
With the notations of the previous example. Then the $i$-th projection map  $T_i:\oplus^n L\rightarrow L, T_i(x_1,\cdots, x_n)=x_i$ is an embedding tensor on $(L, [\c,\c,\c],\alpha)$ with respect to the  representation $(\oplus^n L; \rho, \oplus^n \alpha) $.
\end{example}

\section{The  cohomology of embedding tensors on 3-Hom-Lie algebras}
\def\theequation{\arabic{section}.\arabic{equation}}
\setcounter{equation} {0}

In this section,  we recall some basic results of  representations and cohomologies of
3-Hom-Leibniz algebras.  We construct a representation of the 3-Hom-Leibniz algebra $(V,[\c,\c,\c]_T,\beta)$
on the vector space $L$,  and define the cohomologies of an embedding tensor on 3-Hom-Lie algebras.

\begin{definition}  A representation of the 3-Hom-Leibniz algebra $(\mathcal{L},[\c,\c,\c]_{\mathcal{L}},\alpha)$ is a pair $(V;\beta)$ of
vector space $V$ and a linear map $\beta:V\rightarrow V$, equipped with 3 actions
\begin{align*}
l:\mathcal{L}\otimes \mathcal{L}\otimes V\rightarrow V,\\
m:\mathcal{L}\otimes V\otimes \mathcal{L}\rightarrow V,\\
r:V\otimes \mathcal{L}\otimes \mathcal{L}\rightarrow V,
\end{align*}
satisfying for any  $x,y,z,a,b\in \mathcal{L}$ and $u\in V$
\begin{align*}
l(\alpha(x),\alpha(y),\beta(u))&=\beta(l(x,y,u)),\\
m(\alpha(x),\beta(u),\alpha(z))&=\beta(m(x,u,z)),\\
r(\beta(u),\alpha(y),\alpha(z))&=\beta(r(u,y,z))
\end{align*}
and
\begin{align}
&l([a,b,x]_{\mathcal{L}}, \alpha(y), \beta(u))+ l(\alpha(x),[a,b,y]_{\mathcal{L}},\beta(u))+l(\alpha(x),\alpha(y),l(a,b,u))\nonumber\\
&=l(\alpha(a),\alpha(b),l(x,y,u)),\label{3.1}\\
 & m([a,b,x]_{\mathcal{L}}, \beta(u), \alpha(z))+m(\alpha(x),l(a,b,u),\alpha(z))+m(\alpha(x),\beta(u),[a,b,z]_{\mathcal{L}})\nonumber\\
 &=l(\alpha(a),\alpha(b),m(x,u,z)),\label{3.2}\\
  &r(l(a,b,u), \alpha(y), \alpha(z))+ r(\beta(u),[a,b,y]_{\mathcal{L}},\alpha(z))+r(\beta(u),\alpha(y),[a,b,z]_{\mathcal{L}})\nonumber\\
  &=l(\alpha(a),\alpha(b),r(u,y,z)),\label{3.3}\\
   &r(m(a,u,x), \alpha(y), \alpha(z))+ m(\alpha(x),m(a,u,y),\alpha(z))+l(\alpha(x),\alpha(y),m(a,u,z))\nonumber\\
   &=m(\alpha(a),\beta(u),[x,y,z]_{\mathcal{L}}),\label{3.4}\\
    &r(r(u,b,x), \alpha(y), \alpha(z))+ m(\alpha(x),r(u,b,y),\alpha(z))+l(\alpha(x),\alpha(y),r(u,b,z))\nonumber\\
    &=r(\beta(u),\alpha(b),[x,y,z]_{\mathcal{L}}).\label{3.5}
\end{align}
\end{definition}

An $n$-cochain on a 3-Hom-Leibniz algebra $(\mathcal{L},[\c,\c,\c]_{\mathcal{L}},\alpha)$  with coefficients in a representation  $(V;l,m,r,\beta)$
is a linear map
\begin{align*}
f:\overbrace{\wedge^2\mathcal{L}\otimes\cdots \otimes \wedge^2 \mathcal{L}}^{n-1}\otimes \mathcal{L}\rightarrow V,n\geq 1
\end{align*}
such that $\beta\circ f=f\circ(\tilde{\alpha}^{\otimes n-1}\otimes\alpha)$.
The space generated by $n$-cochains is denoted as $\mathcal{C}^n_{\mathrm{3HLei}}(\mathcal{L},V)$.  The coboundary map $\delta$
from $n$-cochains to
$(n+1)$-cochains, for $X_i=x_i\wedge y_i\in \wedge^2 \mathcal{L}, 1\leq i\leq n$ and $z\in \mathcal{L}$, as
\begin{small}
\begin{align*}
&(\delta f)(X_1,X_2, \cdots, X_n,z)\\
=&\sum_{1\leq j<k\leq n}(-1)^jf(\tilde{\alpha}(X_1),\cdots,\widehat{X_j},\cdots,\tilde{\alpha}(X_{k-1}),\alpha(x_k) \wedge[x_j,y_j,y_k]_{\mathcal{L}}+[x_j,y_j,x_k]_{\mathcal{L}}\wedge  \alpha(y_k),\cdots,\tilde{\alpha}(X_n),\alpha(z))\\
&+\sum_{j=1}^n(-1)^jf(\tilde{\alpha}(X_1),\cdots,\widehat{X_j},\cdots,\tilde{\alpha}(X_{n}),[x_j,y_j,z]_{\mathcal{L}})\\
&+\sum_{j=1}^n(-1)^{j+1}l(\tilde{\alpha}^{n-1}(X_j),f(X_1,\cdots,\widehat{X_j},\cdots,X_{n},z))\\
&+(-1)^{n+1}(m(\alpha^{n-1}(x_n), f(X_1,\cdots,X_{n-1},y_n),\alpha^{n-1}(z))+r(f(X_1,\cdots,X_{n-1},x_n), \alpha^{n-1}(y_n), \alpha^{n-1}(z))).
\end{align*}
\end{small}
It was proved in \cite{Makhlouf} that $\delta^2=0$. Thus,
$(\oplus_{n=1}^{+\infty}\mathcal{C}^n_{\mathrm{3HLei}}(\mathcal{L},V),\delta)$ is a cochain complex.  We denote the set of $n$-cocycles by
$\mathcal{Z}^n_{\mathrm{3HLei}}(\mathcal{L},V)$, the set of $n$-coboundaries by $\mathcal{B}^n_{\mathrm{3HLei}}(\mathcal{L},V)$ and the $n$-th cohomology group by $\mathcal{H}^n_{\mathrm{3HLei}}(\mathcal{L},V)=\mathcal{Z}^n_{\mathrm{3HLei}}(\mathcal{L},V)/\mathcal{B}^n_{\mathrm{3HLei}}(\mathcal{L},V)$.

\begin{lemma}  Let $T:V\rightarrow L$   be an  embedding tensor  on the 3-Hom-Lie algebra $(L, [\c,\c,\c],\alpha)$ with
respect to the representation $(V;\rho,\beta)$.  Define  actions
\begin{align*}
l_T:V\otimes V\otimes L\rightarrow L,m_T:V\otimes L\otimes V\rightarrow L,r_T:L\otimes V\otimes V\rightarrow L,
\end{align*}
by
\begin{align*}
l_T(u,v,x)&=[Tu,Tv,x],\\
m_T(u,x,v)&=[Tu,x,Tv]-T\rho(Tu,x)v,\\
r_T(x,u,v)&=[x,Tu,Tv]-T\rho(x,Tu)v,
\end{align*}
for any $u,v\in V,x\in L.$ Then $(L;l_T,m_T,r_T,\alpha)$ is a representation of the 3-Hom-Leibniz algebra $(V,[\c,\c,\c]_T,\beta)$.
\end{lemma}

{\bf Proof.}  For all $u,v,w,s,t\in V$ and $x\in L$, by Eqs. \eqref{2.2} and \eqref{2.6}, we have
\begin{align*}
l_T(\beta(u),\beta(v),\alpha(x))&=[T\beta(u),T\beta(v),\alpha(x)]\\
&=[\alpha(Tu),\alpha(Tv),\alpha(x)]\\
&=\alpha([Tu,Tv,x])=\alpha (l_T(u,v,x)),\\
m_T(\beta(u),\alpha(x),\beta(v))&=[T\beta(u),\alpha(x),T\beta(v)]-T(\rho(T\beta(u),\alpha(x))\beta(v))\\
&=[\alpha(Tu),\alpha(x),\alpha(Tv)]-T(\rho(\alpha(Tu),\alpha(x))\beta(v))\\
&=\alpha([Tu,x, Tv])-T\beta(\rho(Tu,x)v)\\
&=\alpha([Tu,x, Tv])-\alpha(T\rho(Tu,x)v)\\
&=\alpha (m_T(u,x,v)).
\end{align*}
Similarly, we can show that
$r_T(\alpha(x),\beta(u),\beta(v))=\alpha (r_T(x,u,v))$  holds.

By Eqs. \eqref{2.1}, \eqref{2.7} and \eqref{2.8}, we have
\begin{align*}
&l_T([u,v,s]_{T}, \beta(t), \alpha(x))+ l_T(\beta(s),[u,v,t]_{T},\alpha(x))+l_T(\beta(s),\beta(t),l_T(u,v,x))\\
=&[T[u,v,s]_T , T\beta(t), \alpha(x)]+ [T\beta(s),T[u,v,t]_{T},\alpha(x)]+[T\beta(s),T\beta(t),[Tu,Tv,x]]\\
=&[[Tu,Tv,Ts], \alpha(Tt), \alpha(x)]+ [\alpha(Ts),[Tu,Tv,Tt],\alpha(x)]+[\alpha(Ts),\alpha(Tt),[Tu,Tv,x]]\\
=&[\alpha(Tu),\alpha(Tv),[Ts,Tt,x]]\\
=&[T\beta(u),T\beta(v),[Ts,Tt,x]]\\
=&l_T(\beta(u),\beta(v),l_T(s,t,x)),
\end{align*}
which indicates that Eq. \eqref{3.1} holds.

By Eqs. \eqref{2.1}, \eqref{2.3}, \eqref{2.6},  \eqref{2.7} and \eqref{2.8}, we have
\begin{align*}
 & m_T([u,v,s]_{T},\alpha(x),\beta(t))+m_T(\beta(s),l_T(u,v,x),\beta(t))+m_T(\beta(s),\alpha(x),[u,v,t]_{T}) \\
=&[T[u,v,s]_{T},\alpha(x),T\beta(t)]-T\rho(T[u,v,s]_{T},\alpha(x))\beta(t)+[T\beta(s),[Tu,Tv,x],T\beta(t)]\\
&-T\rho(T\beta(s),[Tu,Tv,x])\beta(t)+[T\beta(s),\alpha(x),T[u,v,t]_{T}]-T\rho(T\beta(s),\alpha(x))[u,v,t]_{T} \\
=&[[Tu,Tv,Ts],\alpha(x),\alpha(Tt)]-T\rho([Tu,Tv,Ts],\alpha(x))\beta(t)+[\alpha(Ts),[Tu,Tv,x],\alpha(Tt)]\\
&-T\rho(\alpha(Ts),[Tu,Tv,x])\beta(t)+[\alpha(Ts),\alpha(x),[Tu,Tv,Tt]]-T\rho(\alpha(Ts),\alpha(x))\rho(Tu,Tv)t \\
=&[\alpha(Tu),\alpha(Tv),[Ts,x,Tt]]-T\rho(\alpha(Tu),\alpha(Tv))\rho(Ts,x)t\\
=&[\alpha(Tu),\alpha(Tv),[Ts,x,Tt]]-[\alpha(Tu),\alpha(Tv),\rho(Ts,x)t]\\
=&[T\beta(u),T\beta(v),[Ts,x,Tt]-\rho(Ts,x)t]\\
=&l_T(u,v,m_T(s,x,t)),\\
&r_T(l_T(u,v,x), \beta(s), \beta(t))+ r_T(\alpha(x),[u,v,s]_{T},\beta(t))+r_T(\alpha(x),\beta(s),[u,v,t]_{T})\\
=&[[Tu,Tv,x], T\beta(s), T\beta(t)]-T\rho([Tu,Tv,x], T\beta(s))\beta(t)+ [\alpha(x),T\rho(Tu,Tv)s,T\beta(t)]\\
&-T\rho(\alpha(x),T\rho(Tu,Tv)s)\beta(t)+[\alpha(x),T\beta(s),T\rho(Tu,Tv)t]-T\rho(\alpha(x),T\beta(s))\rho(Tu,Tv)t\\
=&[[Tu,Tv,x],\alpha(Ts), \alpha(Tt)]-T\rho([Tu,Tv,x], \alpha(Ts))\beta(t)+ [\alpha(x),[Tu,Tv,Ts],\alpha(Tt)]\\
&-T\rho(\alpha(x),[Tu,Tv,Ts])\beta(t)+[\alpha(x),\alpha(Ts),[Tu,Tv,Tt]]-T\rho(\alpha(x),\alpha(Ts))\rho(Tu,Tv)t\\
=&[\alpha(Tu),\alpha(Tv),[x,Ts,Tt]]-T\rho(\alpha(Tu),\alpha(Tv))\rho(x,Ts)t\\
=&[\alpha(Tu),\alpha(Tv),[x,Ts,Tt]]-[\alpha(Tu),\alpha(Tv),T\rho(x,Ts)t]\\
=&[T\beta(u),T\beta(v),[x,Ts,Tt]-T\rho(x,Ts)t]\\
=&l_T(\beta(u),\beta(v),r_T(x,s,t)),
\end{align*}
which imply that Eqs. \eqref{3.2} and \eqref{3.3} hold. Similarly, we can prove that Eqs.\eqref{3.4} and \eqref{3.5} are true. Thus, $(L;l_T,m_T,r_T,\alpha)$ is a representation of the 3-Hom-Leibniz algebra $(V,[\c,\c,\c]_T,\beta)$.
 \hfill $\square$

\medskip

Let  $\delta_T: \mathcal{C}^n_{\mathrm{3HLei}}(V,L)\rightarrow \mathcal{C}^{n+1}_{\mathrm{3HLei}}(V,L),n\geq 1$ be  the coboundary operator of the 3-Hom-Leibniz algebra $(V,[\c,\c,\c]_T,\beta)$
 with coefficients in the representation $(L;l_T,m_T,r_T,\alpha)$. More precisely, for all  $f\in \mathcal{C}^n_{\mathrm{3HLei}}(V,L), V_i=u_i\wedge v_i\in \wedge^2 V, 1\leq i\leq n$ and $w\in V$, we have
\begin{small}
\begin{align*}
&(\delta_T f)(V_1,V_2, \cdots, V_n,w)\\
=&\sum_{1\leq j<k\leq n}(-1)^jf(\tilde{\beta}(V_1),\cdots,\widehat{V_j},\cdots,\tilde{\beta}(V_{k-1}),\beta(u_k) \wedge[u_j,v_j,v_k]_{T}+[u_j,v_j,u_k]_{T}\wedge  \beta(v_k),\cdots,\tilde{\beta}(V_n),\beta(w))\\
&+\sum_{j=1}^n(-1)^jf(\tilde{\beta}(V_1),\cdots,\widehat{V_j},\cdots,\tilde{\beta}(V_{n}),[u_j,v_j,w]_{T})\\
&+\sum_{j=1}^n(-1)^{j+1}l_T(\tilde{\beta}^{n-1}(V_j),f(V_1,\cdots,\widehat{V_j},\cdots,V_{n},w))\\
&+(-1)^{n+1}(m_T(\beta^{n-1}(u_n), f(V_1,\cdots,V_{n-1},v_n),\beta^{n-1}(w))+r_T(f(V_1,\cdots,V_{n-1},u_n), \beta^{n-1}(v_n), \beta^{n-1}(w))).
\end{align*}
\end{small}

In particular, for $f\in \mathcal{C}^1_{\mathrm{3HLei}}(V,L):=\{g\in \mathrm{Hom}(V,L)~|~\alpha\circ g=g\circ \beta\}$ and $u,v,w\in V,$
we have
\begin{align*}
(\delta_T f)(u,v,w)=&-f([u,v,w]_T)+l_T(u,v,f(w))+m_T(u,f(v),w)+r_T(f(u),v,w)\\
=&-f(\rho(Tu,Tv)w)+[Tu,Tv,f(w)]+[Tu,f(v),Tw]-T\rho(Tu,f(v))w\\
&+[f(u),Tv,Tw]-T\rho(f(u),Tv)w
\end{align*}

For any  $(a,b)\in \mathcal{C}^0_{\mathrm{3HLei}}(V,L):=\{(x,y)\in  \wedge^2 L~|~\alpha(x)=x,\alpha(y)=y\}$, we define $\delta_T: \mathcal{C}^0_{\mathrm{3HLei}}(V,L)\rightarrow \mathcal{C}^1_{\mathrm{3HLei}}(V,L), (a, b)\mapsto \wp(a,b)$ by
$$\wp(a,b)v=T\rho(a,b)\beta^{-1}(v)-[a,b,T\beta^{-1}(v)], \forall  v\in V,$$
where $\beta:V\rightarrow V$ is a vector space isomorphism.

\begin{proposition}
Let $T:V\rightarrow L$   be an  embedding tensor  on the regular  3-Hom-Lie algebra $(L, [\c,\c,\c],\alpha)$ with
respect to the regular representation $(V;\rho,\beta)$. Then $\delta_T(\wp(a,b))=0,$  that is the composition  $ \mathcal{C}^0_{\mathrm{3HLei}}(V,L)\stackrel{\delta_{T}}{\longrightarrow} \mathcal{C}^1_{\mathrm{3HLei}}(V,L)\stackrel{\delta_{T}}{\longrightarrow} \mathcal{C}^2_{\mathrm{3HLei}}(V,L)$ is the zero
map.
\end{proposition}

{\bf Proof.}  For any $u,v,w\in V,$ by Eqs. \eqref{2.1},  \eqref{2.2},  \eqref{2.3}, \eqref{2.6} and \eqref{2.7} we have
\begin{align*}
&(\delta_T \wp(a,b))(u,v,w)\\
=&-\wp(a,b)\rho(Tu,Tv)w+[Tu,Tv,\wp(a,b)w]+[Tu,\wp(a,b)v,Tw]-T\rho(Tu,\wp(a,b)v)w\\
&+[\wp(a,b)u,Tv,Tw]-T\rho(\wp(a,b)u,Tv)w\\
=&-T\rho(a,b)\rho(\alpha^{-1}(Tu),\alpha^{-1}(Tv))\beta^{-1}(w)+[a,b,T\rho(\alpha^{-1}(Tu),\alpha^{-1}(Tv))\beta^{-1}(w)]\\
&+[Tu,Tv,T\rho(a,b)\beta^{-1}(w)]-[Tu,Tv,[a,b,T\beta^{-1}(w)]]+[Tu,T\rho(a,b)\beta^{-1}(v),Tw]\\
&-[Tu,[a,b,T\beta^{-1}(v)],Tw]-T\rho(Tu,T\rho(a,b)\beta^{-1}(v))w+T\rho(Tu,[a,b,T\beta^{-1}(v)])w\\
&+[T\rho(a,b)\beta^{-1}(u),Tv,Tw]-[[a,b,T\beta^{-1}(u)],Tv,Tw]-T\rho(T\rho(a,b)\beta^{-1}(u),Tv)w\\
&+T\rho([a,b,T\beta^{-1}(u)],Tv)w\\
=&-T\rho(a,b)\rho(\alpha^{-1}(Tu),\alpha^{-1}(Tv))\beta^{-1}(w)+[a,b,[T\beta^{-1}(u),T\beta^{-1}(v),T\beta^{-1}(w)]]\\
&+T\rho(Tu,Tv)\rho(a,b)\beta^{-1}(w)-[Tu,Tv,[a,b,T\beta^{-1}(w)]]+T\rho(Tu,T\rho(a,b)\beta^{-1}(v))w\\
&-[Tu,[a,b,T\beta^{-1}(v)],Tw]-T\rho(Tu,T\rho(a,b)\beta^{-1}(v))w+T\rho(Tu,[a,b,T\beta^{-1}(v)])w\\
&+T\rho(T\rho(a,b)\beta^{-1}(u),Tv)w-[[a,b,T\beta^{-1}(u)],Tv,Tw]-T\rho(T\rho(a,b)\beta^{-1}(u),Tv)w\\
&+T\rho([a,b,T\beta^{-1}(u)],Tv)w\\
=&0.
\end{align*}
Therefore, $\delta_T(\wp(a,b))=0.$   \hfill $\square$

\medskip

Next we define the cohomology theory of an embedding tensor $T$ on the 3-Hom-Lie algebra $(L, [\c,\c,\c],\alpha)$ with
respect to the representation $(V;\rho,\beta)$. For $n\geq 0$,  define the set of $n$-cochains of $T$ by
$\mathcal{C}^n_T(V,L):=\mathcal{C}^n_{\mathrm{3HLei}}(V,L).$  Then $(\oplus_{n=0}^{\infty}\mathcal{C}^n_T(V,L),\delta_T)$ is a cochain complex.
For $n\geq 1$,  we denote the set of $n$-cocycles by
$\mathcal{Z}^n_T(V,L)$, the set of $n$-coboundaries by $\mathcal{B}^n_T(V,L)$ and the $n$-th cohomology group of the embedding tensor $T$ by $\mathcal{H}^n_T(V,L)=\mathcal{Z}^n_T(V,L)/\mathcal{B}^n_T(V,L)$.

\section{ Deformations of embedding tensors on 3-Hom-Lie algebras}
\def\theequation{\arabic{section}.\arabic{equation}}
\setcounter{equation} {0}

In this section, we discuss linear  deformations of embedding tensors on 3-Hom-Lie algebras.

\begin{definition}
 Let $T:V\rightarrow L$   be an  embedding tensor  on the 3-Hom-Lie algebra $(L, [\c,\c,\c],\alpha)$ with
respect to the representation $(V;\rho,\beta)$, and let  $\mathfrak{I}:V\rightarrow L$  be a linear map.  If  $T_t=T+t\mathfrak{I}$ is a
embedding tensor   on    $(L, [\c,\c,\c],\alpha)$ with
respect to the representation $(V;\rho,\beta)$ for all $t$,
we say that $\mathfrak{I}$ generates a linear deformation of the embedding tensor $T$.

\end{definition}

Suppose that $\mathfrak{I}$ generates a linear deformation of the embedding tensor $T$, then we
have
\begin{align*}
T_t\circ\beta =&\alpha\circ T_t,  \\
[T_tu,T_tv, T_tw ]=& T_t\rho(T_tu,T_tv)w,
\end{align*}
for all $u,v,w\in V.$  This is equivalent to the following conditions
\begin{align}
&\mathfrak{I}\circ\beta=\alpha\circ\mathfrak{I}, \label{4.1}\\
&[Tu,Tv,\mathfrak{I}w]+[Tu,\mathfrak{I}v,Tw]+[\mathfrak{I}u,Tv,Tw]=\mathfrak{I}\rho(Tu,Tv)w +T\rho(Tu,\mathfrak{I}v)w+T\rho(\mathfrak{I}u,Tv)w,\label{4.2}\\
&[Tu,\mathfrak{I}v,\mathfrak{I}w]+[\mathfrak{I}u,Tv,\mathfrak{I}w]+[\mathfrak{I}u,\mathfrak{I}v,Tw]=\mathfrak{I}\rho(Tu,\mathfrak{I}v)w +\mathfrak{I}\rho(\mathfrak{I}u,Tv)w+T\rho(\mathfrak{I}u,\mathfrak{I}v)w,  \label{4.3}\\
&[\mathfrak{I}u,\mathfrak{I}v, \mathfrak{I}w ]= \mathfrak{I}\rho(\mathfrak{I}u,\mathfrak{I}v)w,\label{4.4}
\end{align}
for all $u,v,w\in V.$  Thus, $T_t$
is a linear deformation of $T$ if and only if Eqs. \eqref{4.1}-\eqref{4.4} hold. From Eqs. \eqref{4.1}  and \eqref{4.4} it follows that the
map $\mathfrak{I}$ is  an  embedding tensor  on  the 3-Hom-Lie algebra $(L, [\c,\c,\c],\alpha)$ with
respect to the representation $(V;\rho,\beta)$.

\begin{proposition}
Let $T_t=T+t\mathfrak{I}$  is a linear deformation of an  embedding tensor $T$ on a  3-Hom-Lie algebra $(L, [\c,\c,\c],\alpha)$ with
respect to the representation $(V;\rho,\beta)$.
Then  $\mathfrak{I}\in \mathcal{C}_T^1(V,L)$ is a 1-cocycle of the embedding tensor $T$. Moreover, the 1-cocycle $\mathfrak{I}$ is called the infinitesimal of the linear deformation $T_t$ of $T$.
\end{proposition}

{\bf Proof}\ \  We observe that
Eq.  \eqref{4.2}   implies that $\delta_T\mathfrak{I}=0$. \hfill $\square$

Next, we discuss equivalent linear deformations.
\begin{definition}
 Let $T$  be an  embedding tensor  on the  regular 3-Hom-Lie algebra $(L, [\c,\c,\c],\alpha)$ with
respect to the regular representation $(V;\rho,\beta)$.
Two linear deformations $T^1_t=T+t\mathfrak{I}_1$ and $T^2_t=T+t\mathfrak{I}_2$ are said to be equivalent if   there exist  two elements $a,b\in L$ such that $\alpha(a)=a,\alpha(b)=b$ and  the pair $(\mathrm{Id}_L+t\alpha^{-1}(\mathrm{ad}(a,b)),\mathrm{Id}_V+t\beta^{-1}(\rho(a,b)))$ is a homomorphism from $T^2_t$
to $T^1_t$.
\end{definition}

Let us recall from Definition 2.9 that the pair $(\mathrm{Id}_L+t\alpha^{-1}(\mathrm{ad}(a,b)),\mathrm{Id}_V+t\beta^{-1}(\rho(a,b)))$ is a homomorphism from $T^2_t$
to $T^1_t$ if the following conditions are true:

(1) the map $\mathrm{Id}_L+t\alpha^{-1}(\mathrm{ad}(a,b)):L\rightarrow L,x\mapsto x+t\alpha^{-1}(\mathrm{ad}(a,b)x)$ is a 3-Hom-Lie algebra homomorphism,

(2) $(T+t\mathfrak{I}_1)(u+t\beta^{-1}(\rho(a,b)u))=(\mathrm{Id}_L+t\alpha^{-1}(\mathrm{ad}(a,b)))(Tu+t\mathfrak{I}_2u)$,

(3) $\rho(x,y)u+t\beta^{-1}(\rho(a,b)\rho(x,y)u)=\rho(x+t\alpha^{-1}(\mathrm{ad}(a,b)x),y+t\alpha^{-1}(\mathrm{ad}(a,b)y))(u+t\beta^{-1}(\rho(a,b)u)), \forall x,y\in L, u\in V.$

\begin{theorem}     Let $T$  be an  embedding tensor  on a  regular 3-Hom-Lie algebra $(L, [\c,\c,\c],\alpha)$ with
respect to the  regular representation $(V;\rho,\beta)$. If   two linear deformations $T^1_t=T+t\mathfrak{I}_1$ and $T^2_t=T+t\mathfrak{I}_2$ of $T$ are equivalent, then $\mathfrak{I}_1$ and $\mathfrak{I}_2$ define the same cohomology class in $\mathcal{H}^{1}_T(V,L)$.
\end{theorem}

{\bf Proof}\  \ On comparing the coefficients of   $t^1$ from both sides of the above condition (2), we have
\begin{align*}
\mathfrak{I}_2u-\mathfrak{I}_1u=& T\beta^{-1}(\rho(a,b)u)-\alpha^{-1}([a,b,Tu])\\
=& T\rho(a,b)\beta^{-1}(u)-[a,b,T\beta^{-1}(u)]\\
=& \wp(a,b)u\in  \mathcal{B}^1_T(V,L),
\end{align*}
which implies that  $\mathfrak{I}_2$ and  $\mathfrak{I}_1$ belong to in the same cohomology class in $\mathcal{H}^{1}_T(V,L)$.
\hfill $\square$

\begin{center}
 {\bf ACKNOWLEDGEMENT}
 \end{center}

  The paper is  supported by the NSF of China (No. 12161013), Science and Technology Program of Guizhou Province (No. ZK[2022]031), The Scientific Research Foundation of Guizhou University of Finance and Economics (No.2022KYYB08).

\renewcommand{\refname}{REFERENCES}

\end{document}